\documentclass[final,12pt]{elsarticle}
\usepackage{subcaption}
\usepackage{booktabs}
\usepackage{multirow}
\usepackage{graphicx}
\usepackage{threeparttable}
\usepackage{caption}
\usepackage{amssymb}
\usepackage[a4paper,total={7in,10in}]{geometry} 
\usepackage{color}
\usepackage{makecell}           
\usepackage[fleqn]{amsmath}     
\usepackage[normalem]{ulem}     
\usepackage{lineno}
\usepackage{enumitem}
\usepackage{amsmath,amsfonts,amsthm,bm}
\setenumerate[1]{itemsep=0pt,partopsep=0pt,parsep=\parskip,topsep=5pt}
\setitemize[1]{itemsep=0pt,partopsep=0pt,parsep=\parskip,topsep=5pt}
\setdescription{itemsep=0pt,partopsep=0pt,parsep=\parskip,topsep=5pt}

\newtheorem{Lemma}{Lemma}

\newtheorem{Theorem}{Theorem}
\usepackage[colorlinks,linkcolor=black,citecolor=black,urlcolor=black,]{hyperref}
\setcitestyle{authoryear,semicolon,open={(},close={)}}
\usepackage{siunitx}    
\usepackage{etoolbox} 

\newcommand*\linenomathpatch[1]{%
  \cspreto{#1}{\linenomath}%
  \cspreto{#1*}{\linenomath}%
  \csappto{end#1}{\endlinenomath}%
  \csappto{end#1*}{\endlinenomath}%
}

\linenomathpatch{equation}
\linenomathpatch{gather}
\linenomathpatch{multline}
\linenomathpatch{align}
\linenomathpatch{alignat}
\linenomathpatch{flalign}

\journal{Transportation Research}
\begin{document}



\begin{frontmatter}

\title{Optimal breakpoint selection method for piecewise linear approximation}

\author[inst1]{Shaojun Liu}

\affiliation[inst1]{organization={School of Civil and Environmental Engineering, Nanyang Technological University},
                   country={Singapore}}

\begin{abstract}
Piecewise linearization is a key technique for solving nonlinear problems in transportation network design and other optimization fields, in which generating breakpoints is a fundamental task.
This paper proposes an optimal breakpoint selection method, sequential adjusting method (SAM), to minimize the approximation error between the original function and the piecewise linear function with limited number of pieces, applicable to both convex or concave function.
SAM sequentially adjusts the location of breakpoints based on its adjacent breakpoints, and the optimal positions would be reached after several iterations.
The optimality of the method is proved.
Numerical experiments are conducted on the logarithmic function.
\end{abstract}


\end{frontmatter}


\section{Introduction}
Piecewise linearization is a usefully technique for solving optimization problems involving complex non-linear non-convex components .
\cite{misener2010piecewise,wang2010global} proposed piecewise linearization methods to represent the non-convex components with piecewise linear functions.
\cite{liu2015global,tian2021service,liu2024optimal} applied these methods to reformulate their models into mixed integer linear programming to address the network design problems.
However, in these works, the selection of breakpoints did receive much attention, and they simply employed uniform distributed breakpoints.
\par
Some research contribute to the selection of breakpoints to better approximate the original function.
\cite{birolini2021integrated} applied a differential evolution algorithm to determine the breakpoints, which is problem-related method.
\cite{noruzoliaee2018roads,guo2022shared} determined the breakpoints by selecting from a large set of candidate breakpoints.
This method is time-consuming if high accuracy is desired and is inherently suboptimal.
\cite{xu2021electric} recursively inserts breakpoints to the current optimal position until desired error threshold is reached, which is a greedy method and does not guarantee optimality.
\par
To overcome the research gap, we propose a sequential adjusting method (SAM) for the breakpoint selection of single-value function piecewise linear approximation.
The optimality and convergency are proved in theory for convex and concave functions.
\section{Methodology}
\subsection{Approximation error definition}
To ensure the approximation accuracy of piecewise linear function to the original function, we define two types of approximation errors: minimal maximum absolute error and minimal area difference error.
The absolute error criterion directly minimizes the difference in values between the two functions.
In contrast, the minimal area difference error criterion focuses on minimizes the area of discrepancy between the two functions.
These criteria can be chosen based on demand.
\par
\begin{figure}[!htb]
  \centering
  \includegraphics[width=4in]{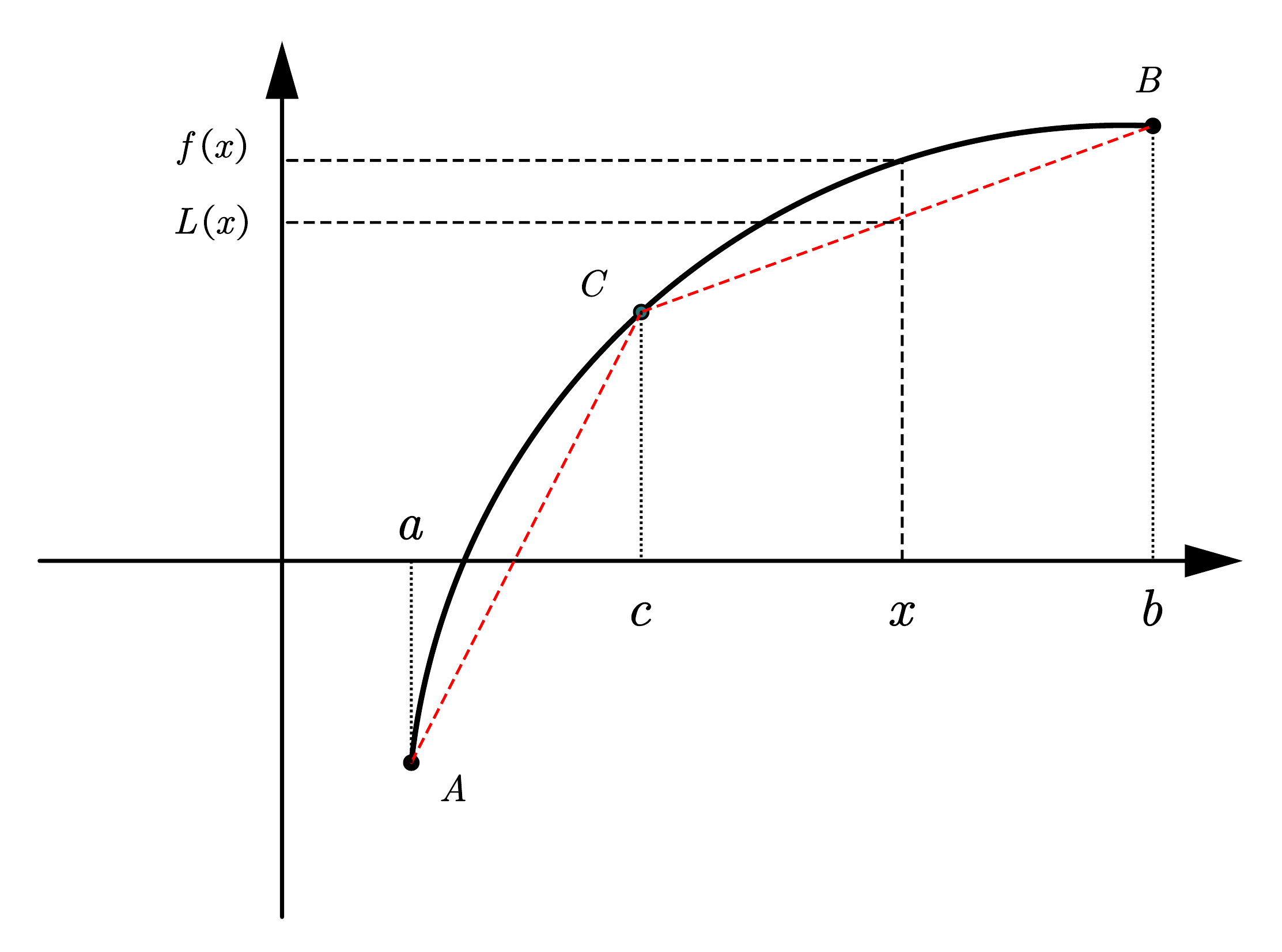}
  \caption{The absolute error}
  \label{fig1.absolute_error}
\end{figure}

The definition of absolute error is shown in Fig.\ref{fig1.absolute_error}.
The absolute error at location $x$ equals to $|f(x)-L(x)|$, where $f(x)$ is the original function value and $L(x)$ is the linear function value.
There exists a $x\in[a,b]$ where the error is maximal, and the minimal maximum absolute error criterion aims to determine the position of the breakpoints (only one point in this example) to minimize this maximum error.
\par
\begin{figure}[!htb]
  \centering
  \includegraphics[width=4in]{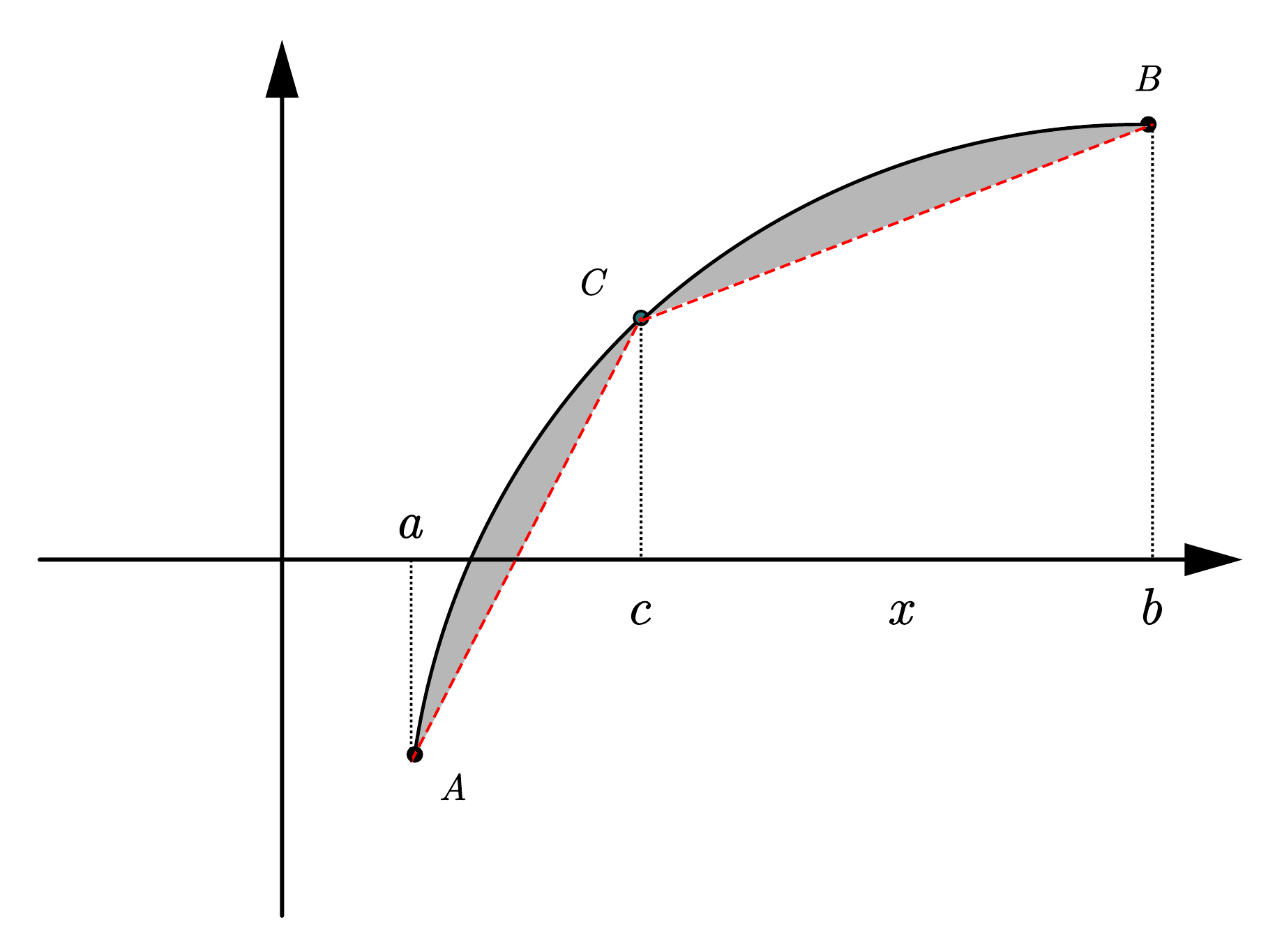}
  \caption{The area difference error}
  \label{fig2.area_difference_error}
\end{figure}

The definition of area difference error is shown in Fig.\ref{fig2.area_difference_error}.
It calculates the total area between the original function and the piecewise linear function (the gray area).
The minimal area difference error criterion aims to determine the position of the breakpoints to minimize the total area.
\par

\subsection{Optimal breakpoints selection for minimal maximum absolute error}
This section discusses the optimal breakpoint selection method with the aim of minimizing the maximum absolute error.
When $f(x)$ is a concave function, $f(x)$ and $-f(x)$ share the same optimal breakpoints due to mirror effect.
In this paper, $f(x)$ is assumed to be concave without any further mention for discussion purposes, and the method is applicable to convex functions.
\subsubsection{The three breakpoints case}
The case with three breakpoints is the simplest scenario, where the positions of the two endpoints are fixed and only the position of the interior point needs to be determined as shown in Fig.\ref{fig1.absolute_error}.
In this case, $a$ and $b$ are known values and our task is to determined the value of $c$.
We formally present the mathematical model for this problem as:
\begin{align}
    \min_{c} \quad&t\\
    s.t. \quad &f(x)-L(x,c)\le t \\
    & x\in[a,b]
\end{align}
which is not a linear programming.
To efficiently obtain the value of $c$, generally one can use numerical methods such as bisection search.
For some special functions, there are analytical solutions.
For example, the optimal value of $c$ for the logarithmic function is $\sqrt{ab}$.
Both the numerical methods and the analytical solution can efficiently find the optimal position.
For convenience of expression, we denote $\varPhi(a,b)$ as the optimal value of $c$ for the minimal maximum absolute error criterion, which can be readily used when $a$ and $b$ are given.

\subsubsection{The general case}
We first present the algorithm of the sequential adjusting method (\textbf{SAM-absolute error}) in Table \ref{table.SAM_absolute}.
The proof will be provided later.
The detailed algorithm steps are presented as below:
\par
\begin{figure}[!htb]
  \centering
  \begin{subfigure}[b]{0.49\textwidth}
  \centering
    \includegraphics[width=0.9\textwidth]{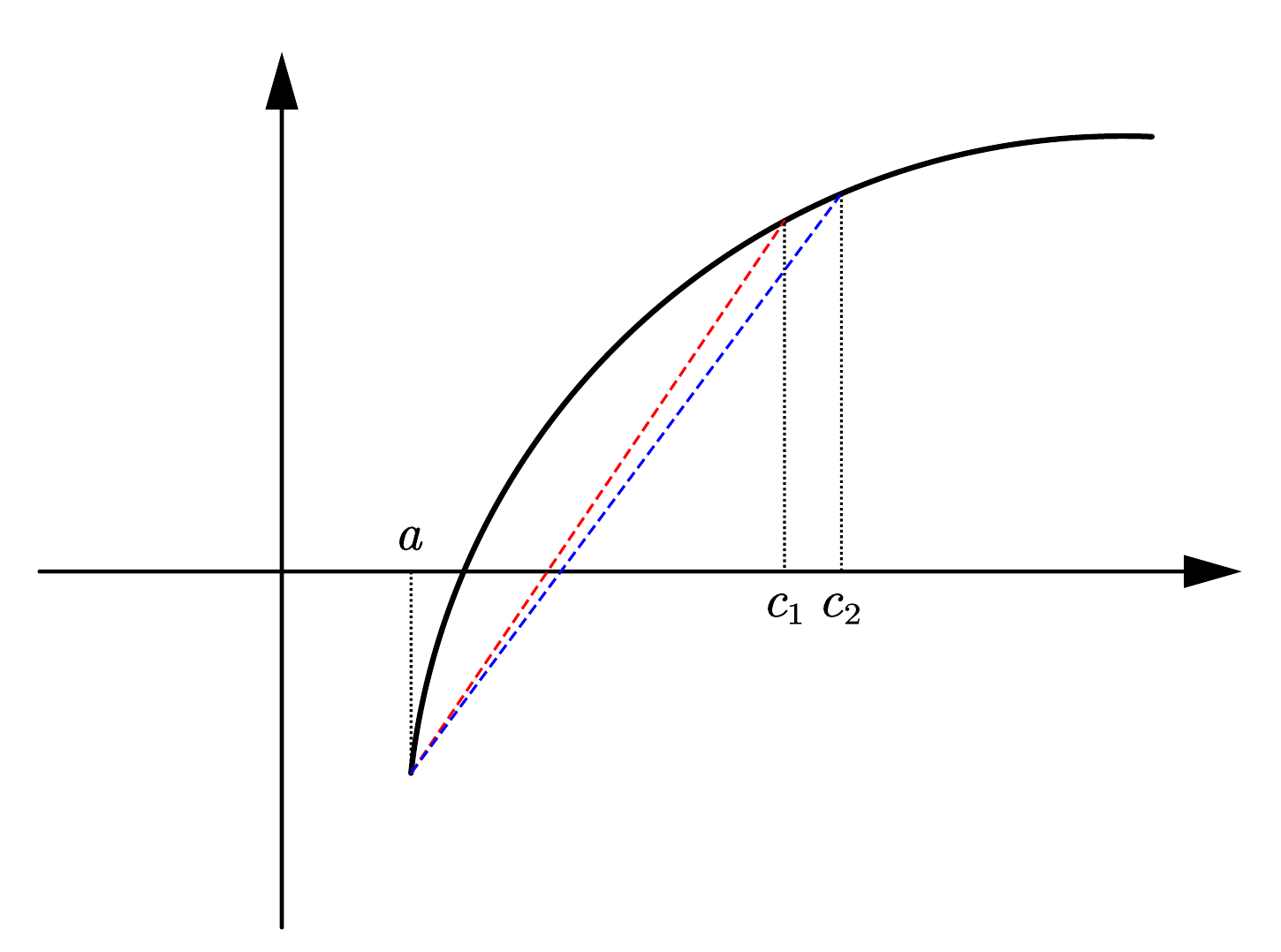}
  \caption{Fixing the left endpoint}
  \label{fig3.1}
  \end{subfigure}
  \hfill
  \begin{subfigure}[b]{0.49\textwidth}
  \centering
    \includegraphics[width=0.9\textwidth]{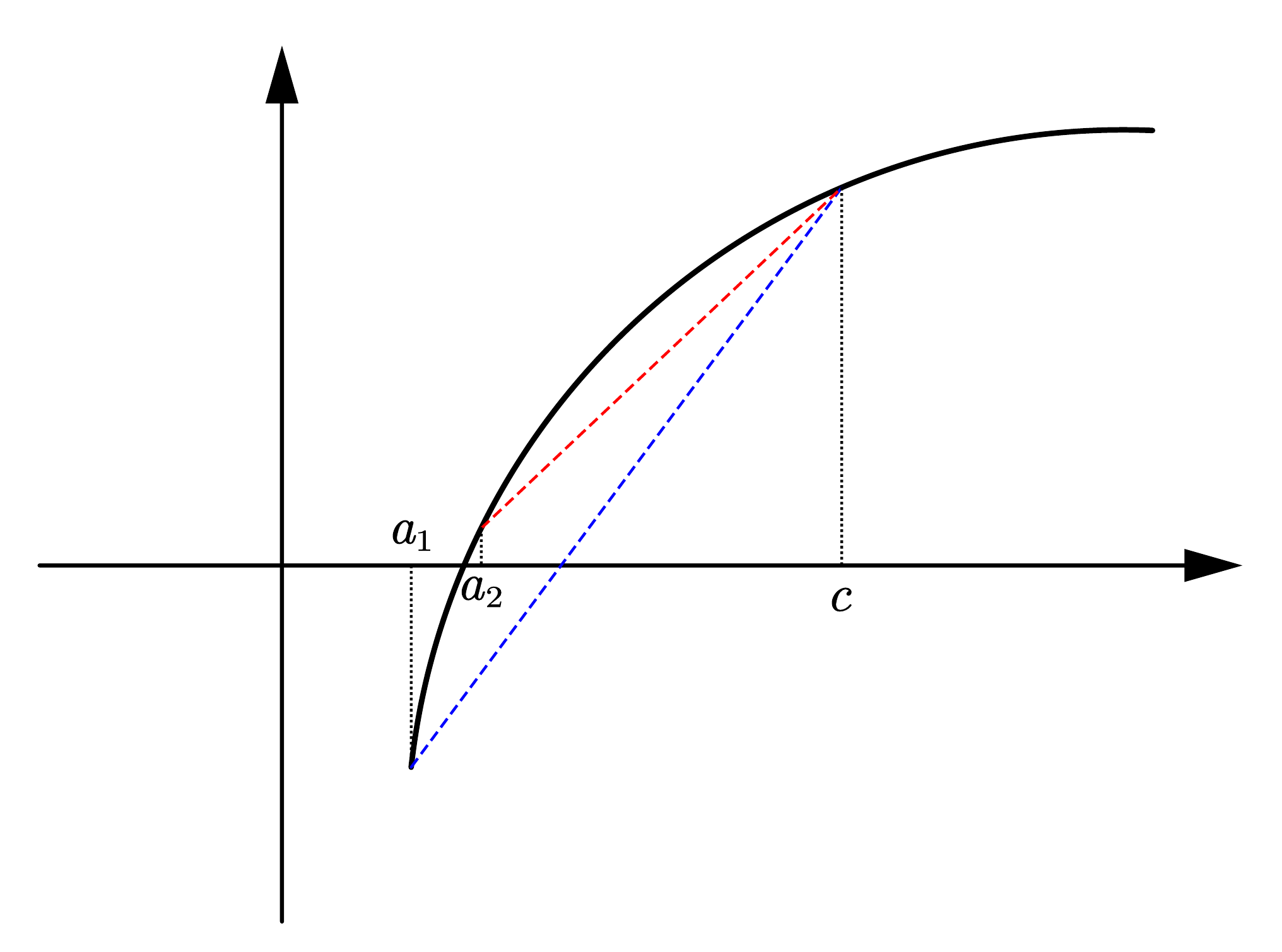}
    \caption{Fixing the right endpoint}
    \label{fig3.2}
  \end{subfigure}
  \caption{Fixing one endpoint of an interval}
  \label{fig:NG}
\end{figure}

\begin{table}[!htb]
  \centering
  \caption{The SAM for the minimal maximum absolute error}
  \label{table.SAM_absolute}
  \resizebox{0.95\textwidth}{!}{
  \begin{tabular}{p{2cm}p{13cm}}
  \toprule
  \textbf{Algorithm:}&{SAM-absolute error}\\
  \midrule
  Step 1: &
  \textbf{\emph{Initialization.}}
  Generate a uniform (or arbitrary) breakpoint set with the pre-given number of breakpoints $N$ and variable range $[\underline{x},\overline{x}]$.
  The set is $bps=(x_1,x_2,...x_N)$, where $x_1=\underline{x},x_N=\overline{x}$ and $x_i<x_{i+1}$.\\
  Step 2: & \textbf{\emph{Updating.}} Update $x_i$ in $bps$ with $\varPhi (x_{i-1},x_{i+1})$ for $i\in (2,3,...,N-1)$ in sequence.
  A new set of breakpoints $bps_{new}$ is then obtained.\\
  Step 3: & \textbf{\emph{Convergence test.}}
  If $abs(bps[i]-bps_{new}[i])\le tolerance$ for all $i\in \{2,3,...,N-1\}$, stop; otherwise, $bps\gets bps_{new}$ and go to Step 2.\\
  \bottomrule
  \end{tabular}  }
\end{table}
The \textbf{SAM-absolute error} algorithm would return the optimal breakpoint set for the minimal maximum absolute error.
To prove it, we present some lemmas first.
\begin{Lemma}\label{lemma_move_end}
  Fixing one endpoint of an interval, the maximal absolute error will increase if the interval length increases and decreases if the interval length decrease.
  The change is also continuous.
  \begin{proof}
   In Fig.\ref{fig3.1}, we fix the left endpoint.
   The red dashed line represents the linear approximation function with endpoint $(a,f(a))$ (simplified as $a$ for clarity, which would be used elsewhere in the absence of ambiguity) and $c_1$ (named as $L(a,c_1)$), while the blue dashed line represents the linear approximation function with endpoint $a$ and $c_2$ (named as $L(a,c_2)$).
   The blue line is to the left of (also below) the red line if $c_1<c_2$, given that $f(x)$ is concave.
   It is clear that the absolute error of $L(a,c_2)$ for all $x\in(a,c_1]$ is large than that of $L(a,c_1)$.
   Therefore, we can conclude that the first part of the lemma is proved for the case where the left endpoint is fixed.
   \par
   Moving $c_1$ to $c_2$ with $\delta$ ($c_2=c_1+\delta$), the continuous changing assumption is equivalent to,
   \begin{align}\label{eq.continous_condtion}
     \lim_{\delta \rightarrow 0} 
     \left\{ 
     \max_{x\in [a,c_2]} \Big( f(x)-L(a,c_2) \Big) 
     -\max_{x\in [a,c_1]} \Big( f(x)-L(a,c_1) \Big) 
     \right\} =0
   \end{align}
   We prove equation \eqref{eq.continous_condtion} by 
   proving that for $\forall x\in [a,c_1]$, the limitation of the difference between $f(x)-L(a,c_1)$ and $f(x)-L(a,c_2)$ is $0$, and that for $\forall x\in [c_1,c_2]$, the limitation of $f(x)-L(a,c_2)$ is $0$.
   \par
   For $x\in [a,c_1]$, we have that
   \begin{align}\label{eq.c1}
     &\Big( f(x)-L(a,c_2) \Big)-\Big( f(x)-L(a,c_1) \Big) \\
     =&L(a,c_1)-L(a,c_2)\\
     =&(x-a)
     \frac
     {[f(c_1)-f(a)](c_1+\delta -a)-[f(c_1+\delta)-f(a)](c_1-a)}
     {(c_1-a)(c_1+\delta-a)}
   \end{align}
   As $\lim_{\delta\to 0}\{[f(c_1)-f(a)](c_1+\delta -a)-[f(c_1+\delta)-f(a)](c_1-a)\}=0$ and 
   $\lim_{\delta\to 0}(c_1-a)(c_1+\delta-a)=(c_1-a)^2$,
   we have $\lim_{\delta\to 0}\left\{\big( f(x)-L(a,c_2) \big)-\big( f(x)-L(a,c_1)\big) \right\}=0$.
   \par
   For $x\in [c_1,c_2]$, we have
   \begin{align}\label{eq.c2}
   &\lim_{\delta\to 0}\{f(x)-L(a,c_2)\}\\
   =& \lim_{\delta\to 0}\{f(x)
   -\frac{f(c_1+\delta)-f(a)}{c_1+\delta-a}(x-c_1-\delta)-f(c_1+\delta)
   \}\\
   =&\lim_{\delta\to 0}\{f(x)
   -\frac{f(c_1)-f(a)}{c_1-a}(x-c_1)-f(c_1)
   \}\\
   =&0
   \end{align}
   \par
   Summarize the above discussion, equation \eqref{eq.continous_condtion} is proved and the maximal absolute error change is continuous when the left endpoint is fixed.
   \par
   The same proof can be conducted when the right endpoint is fixed as shown in Fig.\ref{fig3.2}.
   For brevity, it will not be elaborated further.
   This completed the proof. 
  \end{proof}
\end{Lemma}
\par
\begin{Lemma}\label{lemma_optimal_move}
    Given three points $a< c< b$(as shown in Fig.\ref{fig1.absolute_error}), where $c$ is originally not optimal, 
    let $E_1$ and $E_2$ denote the maximal absolute error for the intervals $[a,c]$ and $[c,b]$, respectively.
    For the optimal point $c^o=\varPhi (a,b)$, denote the interval maximum absolute errors as $E_1^o$ and $E_2^o$.
    It follows that $min(E_1,E_2)<E_1^o=E_2^o < max(E_1,E_2)$.
  \begin{proof}
  We first show that at the optimal point $c^o=\varPhi (a,b)$, $E_1^o\triangleq E_2^o$ holds.
  If at $c=\varPhi (a,b)$, $E_1^o=E_2^o$ does not hold, say $E_1^o<E_2^o$.
  We can move $c$ a enought small distance $\delta$ from $c^o$ to the right, causing $E_1$ to increase, $E_2$ to decrease, and $E_1^o<E_1<E_2<E_2^o$ (by Lemma \ref{lemma_move_end}).
  Consequently, $max(E_1,E_2)$ would decrease, which contradicts the assumption that $\varPhi (a,b)$ is optimal.
  Therefore, $E_1^o=E_2^o$ holds at the optimal point.
  \par
  Since $\varPhi (a,b)$ is optimal and the original $c$ is not optimal, if $c>\varPhi (a,b)$, then we have $E_1>E_1^o$ and $E_2<E_1^o$ according to Lemma \ref{lemma_move_end}.
  Similarly, if $c<\varPhi (a,b)$, then we have $E_1<E_1^o$ and $E_2>E_1^o$.
  Thus,$min(E_1,E_2)<E_1^o=E_2^o< max(E_1,E_2)$ follows.
  This completes the proof.
  \end{proof}
\end{Lemma}

\begin{Theorem}\label{theorem.optimal}
  If $bps=(x_1,x_2,...x_N)$ is the optimal breakpoints to minimize the maximum absolute error, the maximum absolute error for each interval is equal.
  \begin{proof}
  Define $E_i$ as the maximal absolute error for the interval $[x_i, x_{i+1}]$, and $E_{max}$ as the maximal absolute error for $[x_1, x_N]$.
  Then $E_{max}=max\{E_i\}, \forall i \in\{1,2,...N-1\}$.
  \par
  We prove the theorem by contradiction.
  Assume that at the optimum, the breakpoints are $(x_1,x_2,...x_N)$, and not all $E_i$ are equal to $E_{max}$.
  We show that by employing the \textbf{SAM-absolute error} algorithm, $max\{E_i\}$ would be reduced.
  \par
  As not all $E_i$ are equal to $E_{max}$, the first element (indexed as $n$) that makes $E_n=E_{max}$ must lie in one of the following cases: case 1 is $n>1$; case 2 is $n=1$.
  We discuss the two cases separately.
  \par
  Case 1: $n>1$.
  \par
  With \textbf{SAM-absolute error} algorithm, we update $x_i$ in $(x_1,x_2,...x_N)$ with $\varPhi (x_{i-1},x_{i+1})$ for $i\in (2,3,...,N-1)$ in sequence.
  Since $E_i<E_{max}$ (for $i<n$) before updating, and according to Lemma \ref{lemma_optimal_move},
  the updated $E_i$ (for $i<n$) will still be smaller than $E_{max}$.
  As $E_{n-1}<E_{max}$ and $E_{n}=E_{max}$ before updating $x_n$,
  then after updating, $E_{n-1}=E_{n}<E_{max}$ will hold based on  Lemma \ref{lemma_optimal_move}.
  The updating process of $x_n$ tells us that if $E_{i-1}<E_{max}$, after updating, $E_{i-1}$ and $E_i$ would always be smaller than $E_{max}$.
  Repeat this process, all the $E_i$ (for $i>n$) would be smaller than $E_{max}$.
  Summarize the above discuss, the $max\{E_i\}$ would be reduced following the {SAM-absolute error} algorithm if $n>1$.
  \par
  Case 2: $n=1$.
  \par
  $n=1$ means the maximal absolute error of the first interval is $E_{max}$.
  Assume the maximal absolute errors of the first $k$ intervals all equal $E_{max}$, where $k$ could be $1,2,...,N-2$ ($k$ can not be $N-1$ since we have assumed that not all $E_i$ are equal to $E_{max}$).
  \par
  [\textbf{Inner procedure}]: With \textbf{SAM-absolute error} algorithm, when updating the first $k-1$ intervals, $x_i$ does not change because they all have an error of $E_{max}$.
  When it comes to $x_{n+k-1}$ (the last successive interval that has an error of $E_{max}$), since$E_{n+k-1}=E_{max}$ and $E_{n+k}<E_{max}$,
  we get $E_{n+k-1}=E_{n+k}<E_{max}$ after updating.
  The subsequent update falls into Case 1, and the remaining $E_i$ will be smaller than $E_{max}$ after updating.
  \par
  If $k=1$, we have reduced $max\{E_i\}$ with the [\textbf{Inner procedure}].
  Otherwise if $k>1$, we have reduced $k$ by 1 with [\textbf{Inner procedure}]. Repeat the [\textbf{Inner procedure}] $k-1$ times, it would fall into the case of $k=1$.
  \par
  Summarize the above discussion, we have shown that by employing the \textbf{SAM-absolute error} algorithm, $E_{max}$ can be reduced if not all $E_i$ are equal to $E_{max}$.
  At optimality, the maximum absolute error for each interval must be equal.
  This finishes the proof.
\end{proof}
\end{Theorem}
\begin{Theorem}\label{theorem.eqivalent}
If all the interval maximum absolute errors are equal for a solution, $x_i=\varPhi (x_{i-1},x_{i+1})$ would satisfied for all $i\in \{2,3,...,N-1\}$, and vice verse.
\begin{proof}
  We prove the statement that if all the interval maximum absolute errors are equal for a solution, $x_i=\varPhi (x_{i-1},x_{i+1})$ would satisfied for all $i\in \{2,3,...,N-1\}$ by contradiction.
  \par
  If $x_i=\varPhi (x_{i-1},x_{i+1})$ is not satisfied for some $i$, from Lemma \ref{lemma_move_end_area} we get $min\{E_{i-1},E_i\}<max\{E_{i-1},E_i\}$, which implies $E_{i-1}\neq E_i$.
  This proves the first statement.
  \par
  For the statement from the opposite direction, if $x_i=\varPhi (x_{i-1},x_{i+1})$ were satisfied for all $i\in \{2,3,...,N-1\}$, also from Lemma \ref{lemma_move_end_area} we get $E_{i-1}^o=E_{i}^o$.
  Then all the interval maximum absolute errors are equal.
  This completes the proof.
\end{proof}
\end{Theorem}
Now we prove the uniqueness and optimality of the solution whose interval maximum absolute errors are all equal.
\begin{Theorem}\label{theorem.unique}
  The solution in which all the interval maximum absolute errors are equal is both unique and optimal.
  \begin{proof}
  Assume $(x_1,x_2,...,x_{N-1},x_N)$ is a solution in which all the interval maximum absolute errors are equal, with $x_1$ and $x_N$ are fixed endpoints.
  If it is not unique, there must exist another solution $(x_1,y_2,...,y_{N-1},x_N)$, and exist at least one element $x_i\neq y_i$ for $i\in \{2,3,...,N-1\}$.
  Let us assume the index of the first different element is $n$.
  \par
  Case 1: $n>2$.
  In this case, since $x_2=y_2$, then we have $E_1^x=E_1^y$ and it follows that $E_{max}^x$ is equal to $E_{max}^y$ according to Theorem \ref{theorem.optimal}.
  To ensure $E_2^x=E_2^y$, $x_3$ must be equal to $y_3$ ( by Lemma \ref{lemma_move_end}), and this process continues and $x_{i}=y_{i}$ stands for all $i$.
  Consequently, the situation $n>2$ does not occur.
  \par
  Case 2: $n=2$.
  In this case $x_2\neq y_2$.
  Without loss of generality, assume $x_2<y_2$.
  Then we have $E_1^x<E_1^y$ (Lemma \ref{lemma_move_end}) and it follows that $E_{max}^x<E_{max}^y$ (Theorem \ref{theorem.optimal}).
  To ensure $E_2^x<E_2^y$, $x_3$ must be smaller than $y_3$.
  Otherwise, the interval $[y_2,y_3]$ would be a proper subset of $[x_2,x_3]$, and in this situation $E_2^x$ cannot be smaller than $E_2^y$ according to Lemma \ref{lemma_move_end}.
  Continuing this process, we get $x_{N-1}<y_{N-1}$.
  However, this would lead to the final interval $[y_{N-1},x_N]$ being a proper subset of $[x_{N-1},x_N]$, implying $E_{N-1}^x>E_{N-1}^y$ (Lemma \ref{lemma_move_end}), which contradicts $E_{max}^x<E_{max}^y$.
  Consequently, the case $n=2$ does not occur.
  \par
  Summarize the discussion, a second solution that makes all the interval maximum absolute errors equal does not exist, and the uniqueness is therefore proved.
  \par
  As demonstrated by Theorem \ref{theorem.optimal}, the optimal solution must ensure that all its interval maximum absolute errors are equal.
  In the preceding discussion, we proved that the solution satisfying this condition is unique.
  Therefore, a solution in which all the interval maximum absolute errors are equal must optimal.
  This finishes the proof.
  \end{proof}
\end{Theorem}
\par
Theorem \ref{theorem.optimal} has proved that at the optimum, all the maximum absolute error for each interval is equal.
Theorem \ref{theorem.unique} shows that a solution in which all the interval maximum absolute errors are equal is unique and optimal.
Then once we get a solution that all the interval maximum absolute errors are equal, it must be optimal.
The proof of Theorem \ref{theorem.optimal} also demonstrates that the \textbf{SAM-absolute error} algorithm can continually reduce $E_{max}$ if optimality (all $E_i$ are equal) has not been reached.
Since $E_{max}$ must be a positive value, it follows that as the iteration number approaches infinity, the algorithm will eventually find the optimal solution.

\subsection{Optimal breakpoints selection for minimal area difference error}
This section discusses the optimal breakpoint selection method with the aim of minimizing the area difference error.
In this section, we further assume that $f(x)$ is strictly concave.
\subsubsection{Three breakpoints case}
The case with three breakpoints under the minimal area difference error criterion is the scenario, where the positions of the two endpoints are fixed and the position of the interior  point needs to be determined as shown in Fig.\ref{fig2.area_difference_error}.
In this case, $a$ and $b$ are known values and our task is to determined the value of $c$.
\par
\begin{figure}[!htb]
  \centering
  \includegraphics[width=4in]{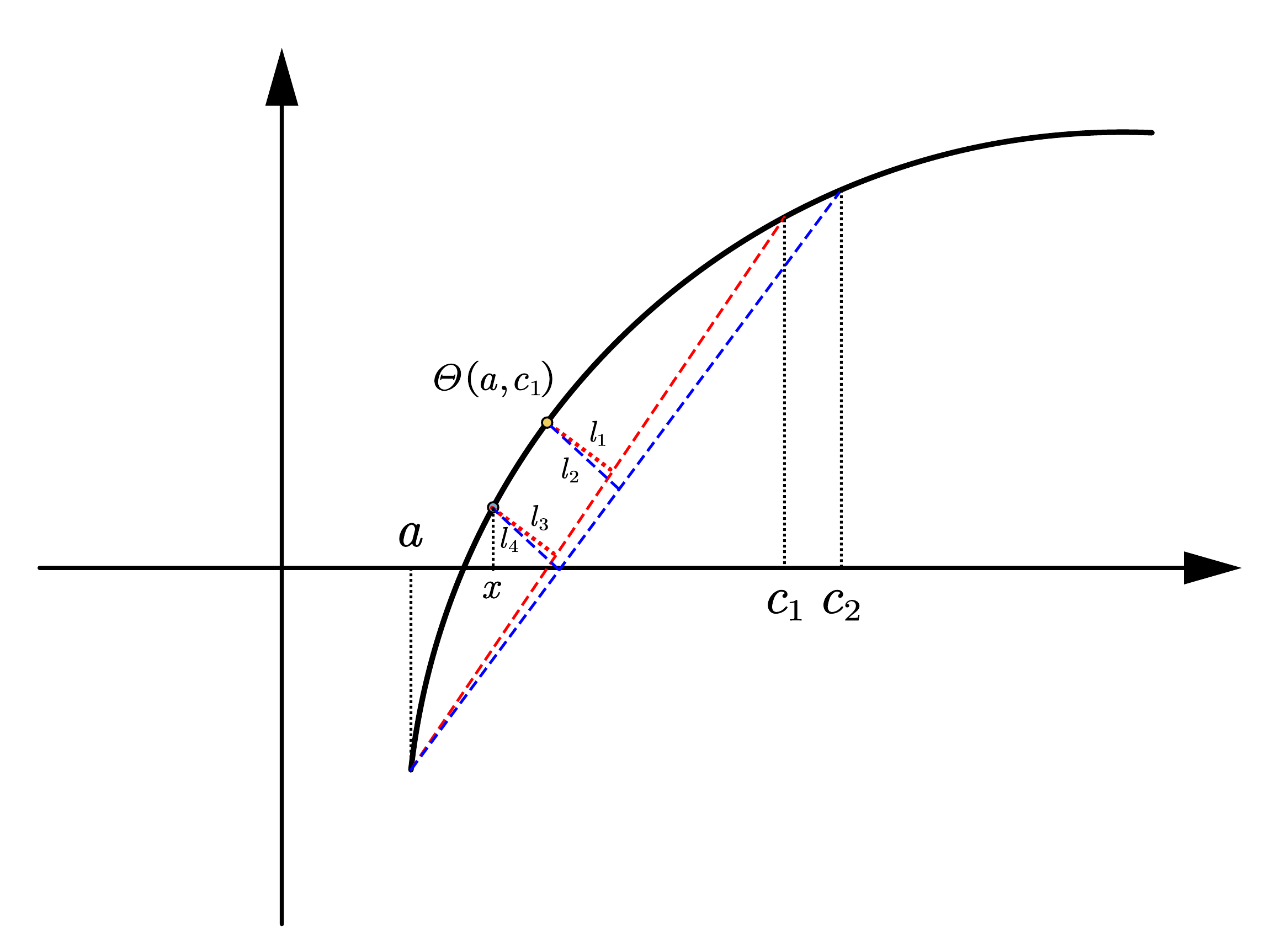}
  \caption{Fixing the left endpoint}
  \label{fig.lemma_move_area}
\end{figure}

Denote $\varTheta(a,b)$ as the optimal value of $c$ under the minimal area difference error criterion.
Lemma \ref{lemma_unique_maxvalue} proves the uniqueness of $\varTheta(a,b)$.
For the general functions,one can employ some numerical methods to find $\varTheta(a,b)$.
For a differentiable function, the unique optimal position should satisfy $f(c)^\prime=\frac{f(b)-f(a)}{b-a}$, and one can find $c$ by methods such as bisection search.
For some special functions, there are analytical solutions.
The optimal position for the logarithmic function is $c=\frac{b-a}{f(b)-f(a)}$.
For convenience of expression, we assume that $\varTheta(a,b)$ can be readily used when $a$ and $b$ are given.
\begin{Lemma}\label{lemma_unique_maxvalue}
    Given three points $a< c< b$(as shown in Fig.\ref{fig2.area_difference_error}, the value of $\varTheta(a,b)$ is unique if $f(x)$ is strictly concave.
  \begin{proof}
  As minimizing the area of the gray zone is equivalent to maximizing the area of the triangle $ABC$, since the total area between the original function and the linear function is fixed.
  To maximize the area of the triangle $ABC$, we consider a function $g(x)=f(x)-L(a,b),x\in [a,b]$, where $L(a,b)$ is the linear function across $(a,f(a))$ and $(b,f(b))$.
  The area of $ABC$ equals $\frac{g(x)\cdot |AB|}{2}$.
  As $g(x)$ is a strict concave function if $f(x)$ is strictly concave, its maximum value is unique.
  Consequently, the value of $\varTheta(a,b)$ is unique.
  \end{proof}
\end{Lemma}
\begin{Lemma}\label{lemma_move_end_area}
    For intervals $[a,c_1]$ and $[a,c_2]$, if $c_1<c_2$, then $\varTheta (a,c_1)\le \varTheta (a,c_2)$;
    for intervals $[a_1,c]$ and $[a_2,c]$, if $a_1<a_2$, then $\varTheta (a_1,c)\le \varTheta (a_2,c)$.
  \begin{proof}
  We first prove that if $c_1<c_2$, then $\varTheta (a,c_1)\le \varTheta (a,c_2)$.
  As shown in Fig.\ref{fig.lemma_move_area}, assume $l_1$ and $l_2$ are the distances from the optimal position $\varTheta (a,c_1)$ to the two linear functions $L(a,c_1)$ and $L(a,c_2)$, respectively.
  For any $x<\varTheta (a,c_1)$, let $l_3$ and $l_4$ be the distances to the two linear functions.
  We know $l_1$ is parallel to $l_2$, and $l_3$ is parallel to $l_4$.
  Since $l_3<l_1$ (by Lemma \ref{lemma_unique_maxvalue}), it follows that $l_4<l_2$.
  The distance between $\varTheta (a,c_2)$ and $L(a,c_2)$ (which defined as the largest distance) cannot be smaller than $l_2$.
  Therefore, if $x<\varTheta (a,c_1)$, $x$ cannot be $\varTheta (a,c_2)$, which implies $\varTheta (a,c_1)\le \varTheta (a,c_2)$.
  \par
  The case where $a_1<a_2$ and $\varTheta (a_1,c)\le \varTheta (a_2,c)$ can be proved using the same method.
  This completes the proof.
  \end{proof}
\end{Lemma}
\subsubsection{The general case}
We present the algorithm of the sequential adjusting method (\textbf{SAM-area error}) for minimal area difference error in Table \ref{table.SAM_area}.
The proof will be provided later.
The detailed algorithm steps are presented as below:
\par
\begin{table}[!htb]
  \centering
  \caption{The SAM for the minimal area difference error}
  \label{table.SAM_area}
  \resizebox{0.95\textwidth}{!}{
  \begin{tabular}{p{2cm}p{13cm}}
  \toprule
  \textbf{Algorithm:}&{SAM-area error}\\
  \midrule
  Step 1: &
  \textbf{\emph{Initialization.}}
  Generate a uniform (or arbitrary) breakpoint set with the pre-given number of breakpoints $N$ and variable range $[\underline{x},\overline{x}]$.
  The set is $bps=(x_1,x_2,...x_N)$, where $x_1=\underline{x},x_N=\overline{x}$ and $x_i<x_{i+1}$.\\
  Step 2: & \textbf{\emph{Updating.}} Update $x_i$ in $bps$ with $\varTheta (x_{i-1},x_{i+1})$ for $i\in (2,3,...,N-1)$ in sequence.
  A new set of breakpoints $bps_{new}$ is then obtained.\\
  Step 3: & \textbf{\emph{Convergence test.}}
  If $abs(bps[i]-bps_{new}[i])\le tolerance$ for all $i\in \{2,3,...,N-1\}$, stop; otherwise, $bps\gets bps_{new}$ and go to Step 2.\\
  \bottomrule
  \end{tabular}  }
\end{table}
The \textbf{SAM-area error} algorithm would return an optimal breakpoint set for the minimal area difference error.
\begin{Theorem}\label{theorem.optimality_area}
  If $bps=(x_1,x_2,...x_N)$ is the optimal breakpoints to minimize the area difference error, $x_i=\varTheta (x_{i-1},x_{i+1})$ must be satisfied for all $i\in \{2,3,...,N-1\}$.
  \begin{proof}
  Assume $(x_1,x_2,...,x_{N-1},x_N)$ is a optimal solution, where $x_1,x_N$ are fixed endpoints.
  If $x_i \neq \varTheta  (x_{i-1},x_{i+1})$ for some $i\in \{2,3,...,N-1\}$.
  Then by moving $x_i$ to $\varTheta  (x_{i-1},x_{i+1})$, the total area difference would be reduced.
  Therefore,$x_i=\varTheta  (x_{i-1},x_{i+1})$ must be satisfied for all $i\in \{2,3,...,N-1\}$.
  This completes the proof.
  \end{proof}
\end{Theorem}
\begin{Theorem}\label{theorem.unique_area}
  The breakpoints $(x_1,x_2,...,x_{N-1},x_N)$ that satisfying $x_i=\varPhi (x_{i-1},x_{i+1})$ for all $i\in \{2,3,...,N-1\}$ is unique and optimal.
  \begin{proof}
  Assume $(x_1,x_2,...,x_{N-1},x_N)$ satisfies $x_i=\varTheta  (x_{i-1},x_{i+1})$ for all $i\in \{2,3,...,N-1\}$ , where $x_1,x_N$ are fixed endpoints.
  If such breakpoints is not unique, there must exist another solution $(x_1,y_2,...,y_{N-1},x_N)$, and exist at least  one element $x_i\neq y_i$ for $i\in \{2,3,...,N-1\}$.
  Let the index of the first different element is $n$.
  \par
  Without loss of generality, assume $x_n<y_n$.
  According to Lemma \ref{lemma_unique_maxvalue} and Lemma \ref{lemma_move_end_area}, $x_{n+1}$ cannot be equal to or larger than $y_{n+1}$
  Therefore we have $x_{n+1}<y_{n+1}$.
  This process continues, leading to $x_N<y_N$, where $x_N$ and $y_N$ are endpoints and should be equal.
  \par
  Therefore, a second solution does not exist, and the uniqueness is proved.
  \par
  As demonstrated by Theorem \ref{theorem.optimality_area}, the optimal solution must satisfy that $x_i=\varTheta (x_{i-1},x_{i+1})$ holds for all $i\in \{2,3,...,N-1\}$.
  The above discussion has proved such solution is unique.
  Therefore, the solution that ensures $x_i=\varPhi (x_{i-1},x_{i+1})$ for all $i\in \{2,3,...,N-1\}$ must be optimal.
  This completes the proof.
  \end{proof}
\end{Theorem}
Theorem \ref{theorem.optimality_area} and Theorem \ref{theorem.unique_area} have proved that, at optimality, $x_i=\varTheta(x_{i-1},x_{i+1})$ for all $i$ and the optimal solution is unique.
The \textbf{SAM-area error} algorithm would always reduce the area difference error by updating $x_i\gets  \varTheta(x_{i-1},x_{i+1})$ when they are not equal.
Since the area difference error must be a positive value, it follows that as the iteration number approaches infinity, the \textbf{SAM-area error} algorithm will eventually reach the optimal solution.

\section{Numerical experiments}
In this section, we test the two algorithms: \textbf{SAM-absolute error} and \textbf{SAM-area error} on the logarithmic function and show the convergency.
The given feasible range is $[0.1,10]$ and the desired number of breakpoints is $5$.
The tolerance for convergency checking is $10^{-8}$.
We would generate a uniform breakpoints set as the initial set, and then employ the two algorithms to find the optimal solution.
For the logarithmic function with an interval $[a,b]$, the optimal location to place the interior point is $\varPhi(a,b)=\sqrt{ab}$ under the minimal maximum absolute error criterion, and $\varTheta (a,b) = \frac{b-a}{ln(b)-ln(a)}$ under the minimal area difference error criterion.
We would employ the two algorithms and calculate both the maximal absolute error and area difference error for comparison, even though one is not the optimization objective for each case.
\par
\begin{figure}[!htb]
  \centering
  \includegraphics[width=3.5in]{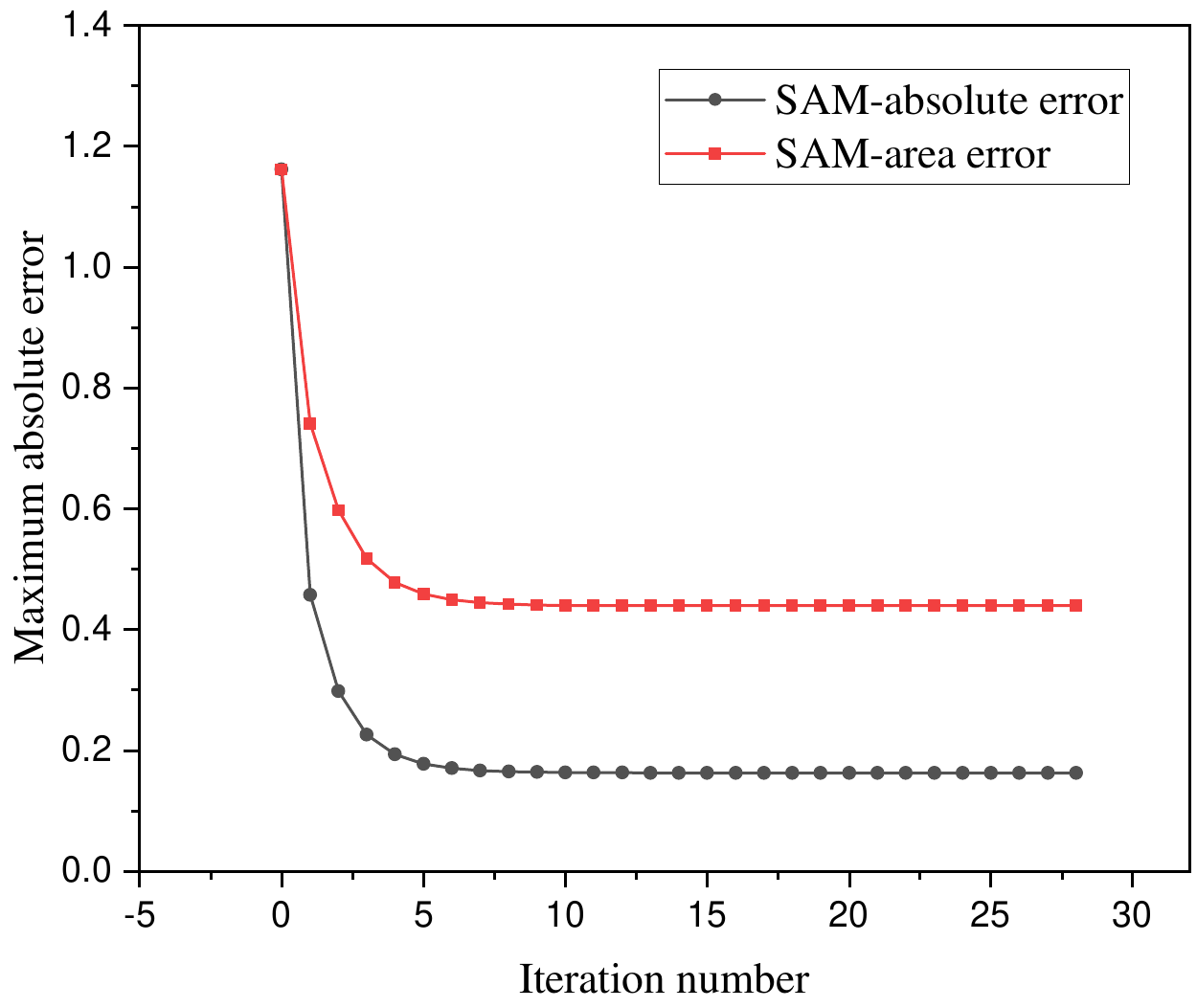}
  \caption{Maximum absolute error}
  \label{fig.e_max}
\end{figure}

Fig.\ref{fig.e_max} displays the convergence process of the maximum absolute error for the two algorithms.
In the figure, the black line represents the results of \textbf{SAM-absolute error} algorithm, which consistently decreases before converging and is smaller than that of the \textbf{SAM-area error} algorithm after convergence.
\par
\par
\begin{figure}[!htb]
  \centering
  \includegraphics[width=3.5in]{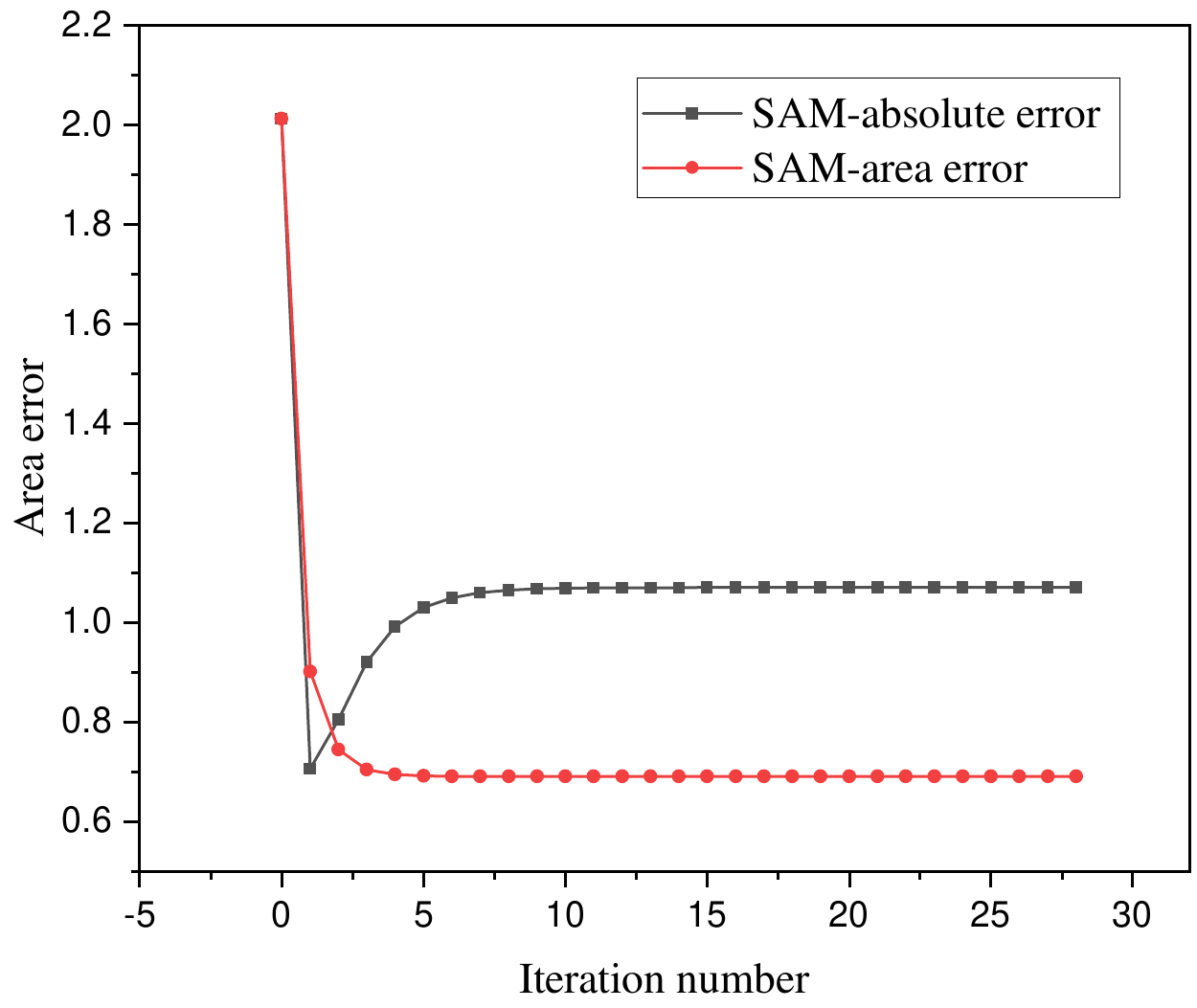}
  \caption{Area error}
  \label{fig.area}
\end{figure}

Fig.\ref{fig.area} shows the convergence of area error.
In the figure, the red line represents the results of \textbf{SAM-area error} algorithm, which consistently decreases before converging and is smaller than that of the \textbf{SAM-absolute error} algorithm after convergence.
It can be found that with the \textbf{SAM-absolute error} algorithm, the area error could increase in the iteration process, since it is the optimization objective of this algorithm.
\par
\par
Combining the results in Fig.\ref{fig.e_max} and Fig.\ref{fig.area}, the \textbf{SAM-absolute error} algorithm achieves a better maximum absolute error solution, while the \textbf{SAM-area error} algorithm provides a superior area error solution.
These two objectives generally conflict with each other.
\par
\begin{figure}[!htb]
  \centering
  \includegraphics[width=3.5in]{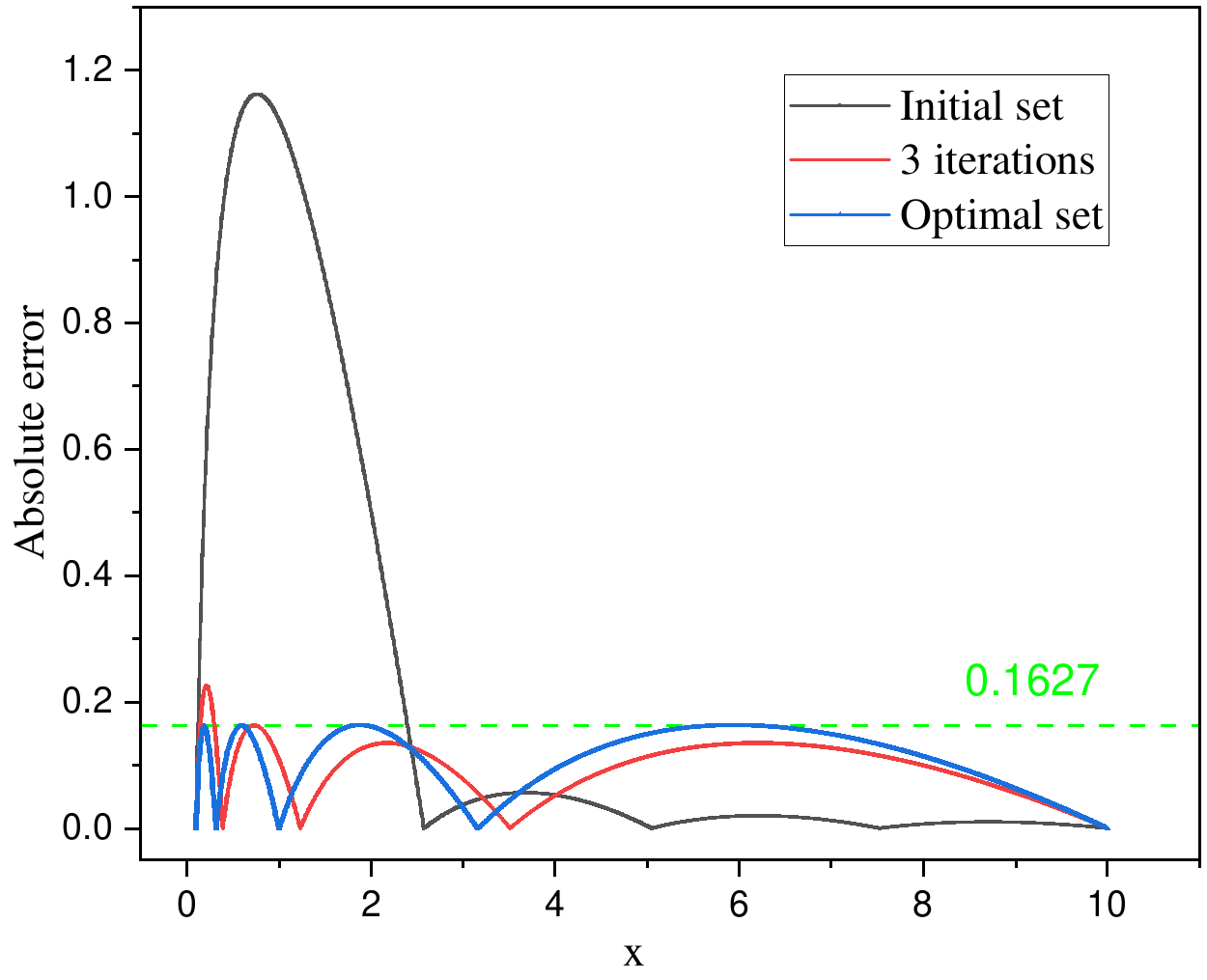}
  \caption{Approximate error of \textbf{SAM-absolute error} algorithm}
  \label{fig.ram1.error}
\end{figure}

Fig.\ref{fig.ram1.error} presents the approximation errors (calculated by $ln(x)-L_i$) for the \textbf{SAM-absolute error} algorithm with different breakpoint sets: the black line represents the initial uniform breakpoints, the red line indicates  the breakpoints after three iterations, and the blue line is the optimal breakpoints.
It can be seen that the maximum absolute error is significantly  reduced from $1.16$ to $0.25$ after only three iterations, and it is close to the optimal value $0.16$.
At optimality, all the maximum absolute error for different intervals are equal, which consists with Theorem \ref{theorem.optimal}.
\par
\par
\begin{figure}[!htb]
  \centering
  \includegraphics[width=3.5in]{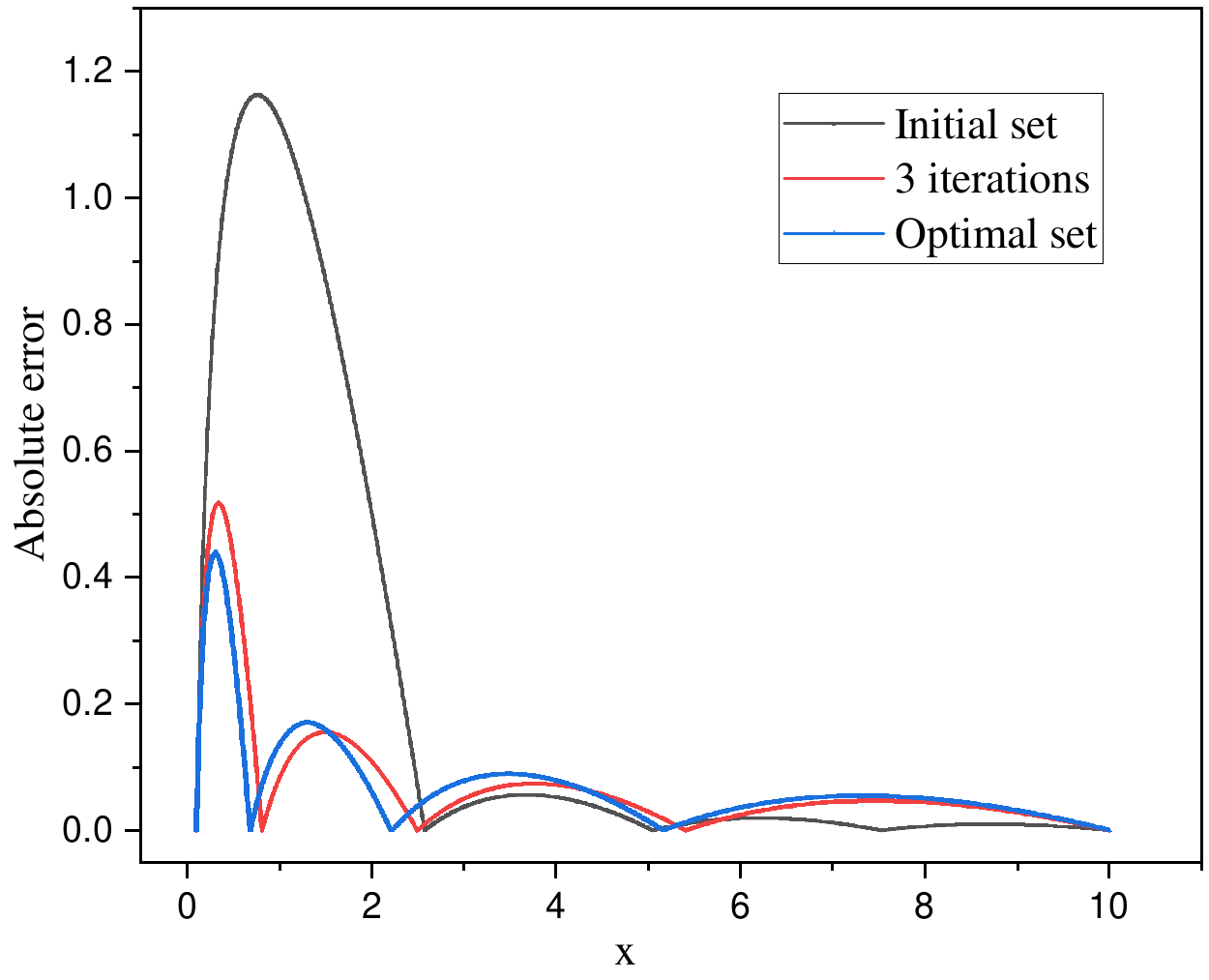}
  \caption{Approximate error of \textbf{SAM-area error} algorithm}
  \label{fig.ram2.error}
\end{figure}

Fig.\ref{fig.ram2.error} illustrates the approximation errors for the \textbf{SAM-area error} algorithm with different breakpoint sets, where the black line denotes the initial uniform breakpoints, the red line represents the breakpoints after three iterations, and the blue line is the optimal breakpoints.
Comparing with the optimal solution of \textbf{SAM-absolute error} algorithm in Fig.\ref{fig.ram1.error}, the maximum absolute error for different intervals are generally not equal.
\par
To demonstrate the efficiency of the algorithms, we evaluate algorithms with different number of breakpoints, ranging from $5$ to $100$, and the resulting running time is displayed in Table \ref{table.time_test}.
Only the running time results of the \textbf{SAM-absolute error} algorithm are presented, since the results of \textbf{SAM-area error} are quite close.
The results indicate that the running time is very short.
Given the nature that only a limited number of breakpoints would be adopted in piecewise linearization methods, this \textbf{SAM} approach can be employed without concern.
\par
\begin{table}[!htb]
\centering
\caption{Running time for different breakpoint number}\label{table.time_test}
\begin{tabular}{ccccccccc}
    \hline $N$	& $5$&$10$&$20$&$50$&$100$\\
    \hline
    time($s$)    & $5.16\times10^{-4}$& $4.48\times10^{-3}$ &$3.28\times10^{-2}$&$0.49$&$3.63$\\
    \hline
\end{tabular}
\end{table}
\section{Conclusions}
This paper proposes SAM methods for breakpoint selection in piecewise linear approximation.
The optimality, uniqueness and convergence of the methods are proved.
Numerical experiments demonstrate that the algorithm is of high efficiency, and optimal results can be achieved with a limited number of iterations.
Although SAM is proved of optimality only for convex functions and concave functions, it could also be applied on non-convex or non-concave functions since they could also get a good solution in experience.
\section*{Acknowledgments}
This paper is supported by the Singapore MOE AcRF Tier 1 project MOE2021-T1-002-062.

\bibliographystyle{elsarticle-harv} 

\bibliography{elsarticle_brue}

\begin{thebibliography}{9}
\expandafter\ifx\csname natexlab\endcsname\relax\def\natexlab#1{#1}\fi
\providecommand{\url}[1]{\texttt{#1}}
\providecommand{\href}[2]{#2}
\providecommand{\path}[1]{#1}
\providecommand{\DOIprefix}{doi:}
\providecommand{\ArXivprefix}{arXiv:}
\providecommand{\URLprefix}{URL: }
\providecommand{\Pubmedprefix}{pmid:}
\providecommand{\doi}[1]{\href{http://dx.doi.org/#1}{\path{#1}}}
\providecommand{\Pubmed}[1]{\href{pmid:#1}{\path{#1}}}
\providecommand{\bibinfo}[2]{#2}
\ifx\xfnm\relax \def\xfnm[#1]{\unskip,\space#1}\fi
\bibitem[{Birolini et~al.(2021)Birolini, Antunes, Cattaneo, Malighetti and
  Paleari}]{birolini2021integrated}
\bibinfo{author}{Birolini, S.}, \bibinfo{author}{Antunes, A.P.},
  \bibinfo{author}{Cattaneo, M.}, \bibinfo{author}{Malighetti, P.},
  \bibinfo{author}{Paleari, S.}, \bibinfo{year}{2021}.
\newblock \bibinfo{title}{Integrated flight scheduling and fleet assignment
  with improved supply-demand interactions}.
\newblock \bibinfo{journal}{Transportation Research Part B: Methodological}
  \bibinfo{volume}{149}, \bibinfo{pages}{162--180}.
\bibitem[{Guo et~al.(2022)Guo, Chen and Liu}]{guo2022shared}
\bibinfo{author}{Guo, H.}, \bibinfo{author}{Chen, Y.}, \bibinfo{author}{Liu,
  Y.}, \bibinfo{year}{2022}.
\newblock \bibinfo{title}{Shared autonomous vehicle management considering
  competition with human-driven private vehicles}.
\newblock \bibinfo{journal}{Transportation research part C: emerging
  technologies} \bibinfo{volume}{136}, \bibinfo{pages}{103547}.
\bibitem[{Liu and Wang(2015)}]{liu2015global}
\bibinfo{author}{Liu, H.}, \bibinfo{author}{Wang, D.Z.}, \bibinfo{year}{2015}.
\newblock \bibinfo{title}{Global optimization method for network design problem
  with stochastic user equilibrium}.
\newblock \bibinfo{journal}{Transportation Research Part B: Methodological}
  \bibinfo{volume}{72}, \bibinfo{pages}{20--39}.
\bibitem[{Liu et~al.(2024)Liu, Wang, Tian and Lin}]{liu2024optimal}
\bibinfo{author}{Liu, S.}, \bibinfo{author}{Wang, D.Z.}, \bibinfo{author}{Tian,
  Q.}, \bibinfo{author}{Lin, Y.H.}, \bibinfo{year}{2024}.
\newblock \bibinfo{title}{Optimal configuration of dynamic wireless charging
  facilities considering electric vehicle battery capacity}.
\newblock \bibinfo{journal}{Transportation Research Part E: Logistics and
  Transportation Review} \bibinfo{volume}{181}, \bibinfo{pages}{103376}.
\bibitem[{Misener and Floudas(2010)}]{misener2010piecewise}
\bibinfo{author}{Misener, R.}, \bibinfo{author}{Floudas, C.},
  \bibinfo{year}{2010}.
\newblock \bibinfo{title}{Piecewise-linear approximations of multidimensional
  functions}.
\newblock \bibinfo{journal}{Journal of optimization theory and applications}
  \bibinfo{volume}{145}, \bibinfo{pages}{120--147}.
\bibitem[{Noruzoliaee et~al.(2018)Noruzoliaee, Zou and
  Liu}]{noruzoliaee2018roads}
\bibinfo{author}{Noruzoliaee, M.}, \bibinfo{author}{Zou, B.},
  \bibinfo{author}{Liu, Y.}, \bibinfo{year}{2018}.
\newblock \bibinfo{title}{Roads in transition: Integrated modeling of a
  manufacturer-traveler-infrastructure system in a mixed autonomous/human
  driving environment}.
\newblock \bibinfo{journal}{Transportation Research Part C: Emerging
  Technologies} \bibinfo{volume}{90}, \bibinfo{pages}{307--333}.
\bibitem[{Tian et~al.(2021)Tian, Wang and Lin}]{tian2021service}
\bibinfo{author}{Tian, Q.}, \bibinfo{author}{Wang, D.Z.}, \bibinfo{author}{Lin,
  Y.H.}, \bibinfo{year}{2021}.
\newblock \bibinfo{title}{Service operation design in a transit network with
  congested common lines}.
\newblock \bibinfo{journal}{Transportation Research Part B: Methodological}
  \bibinfo{volume}{144}, \bibinfo{pages}{81--102}.
\bibitem[{Wang and Lo(2010)}]{wang2010global}
\bibinfo{author}{Wang, D.Z.}, \bibinfo{author}{Lo, H.K.}, \bibinfo{year}{2010}.
\newblock \bibinfo{title}{Global optimum of the linearized network design
  problem with equilibrium flows}.
\newblock \bibinfo{journal}{Transportation Research Part B: Methodological}
  \bibinfo{volume}{44}, \bibinfo{pages}{482--492}.
\bibitem[{Xu et~al.(2021)Xu, Wu and Tan}]{xu2021electric}
\bibinfo{author}{Xu, M.}, \bibinfo{author}{Wu, T.}, \bibinfo{author}{Tan, Z.},
  \bibinfo{year}{2021}.
\newblock \bibinfo{title}{Electric vehicle fleet size for carsharing services
  considering on-demand charging strategy and battery degradation}.
\newblock \bibinfo{journal}{Transportation Research Part C: Emerging
  Technologies} \bibinfo{volume}{127}, \bibinfo{pages}{103146}.

\end{thebibliography}

\end{document}